\documentclass[12pt]{article}\usepackage{amsmath}\usepackage{amssymb}
 
\author{
Krzysztof Ciesielski%
\thanks{
The work of the first author was partially supported by 
NSF Cooperative
Research Grant INT-9600548, with its Polish part 
financed by Polish Academy of Science PAN.}
\\
{\footnotesize Department of Mathematics,}
{\footnotesize West Virginia University,}\\
{\footnotesize Morgantown, WV 26506-6310, USA}\\
{\footnotesize KCies@wvnvms.wvnet.edu}
\and
Saharon Shelah%
\thanks{This work was supported in part by a grant from
``Basic Research Foundation'' 
founded by the Israel Academy of Sciences and Humanities.
Publication 695. 
}\\
{\footnotesize Institute of Mathematics,}
{\footnotesize the Hebrew University of Jerusalem}\\
{\footnotesize 91904 Jerusalem, Israel}\\
{\footnotesize and}\\
{\footnotesize Department of Mathematics,}
{\footnotesize Rutgers University}\\
{\footnotesize New Brunswick, NJ 08854, USA}}

\title{Category analog of sup-measurability problem}

\date{}

\newcommand{\cont}{{\mathfrak c}}

\newcommand{\forces}{\mathrel{\|}\joinrel\mathrel{-}}

\newcommand{\cf}{{\rm cf}}

\newcommand{\real}{{\mathbb R}}
\newcommand{\R}{{\real}}
\newcommand{\rational}{{\mathbb Q}}

\newcommand{\la}{{\langle}}
\newcommand{\ra}{{\rangle}}

\newcommand{\M}{{\cal M}}

\newcommand{\K}{{\cal K}}

\newcommand{\e}{{\emptyset}}
\newcommand{\ep}{{\varepsilon}}

\newcommand{\proj}{{\rm pr}}

\def\cl{\mathop{\rm cl}}

\def\proof{\noindent {\sc Proof. }}
\def\qed{\hfill\vrule height6pt width6pt depth1pt\medskip}

\renewcommand{\P}{{\cal P}}

\newtheorem{theorem}{Theorem}
\newtheorem{corollary}[theorem]{Corollary}
\newtheorem{proposition}[theorem]{Proposition}
\newtheorem{lemma}[theorem]{Lemma}
\newtheorem{problem}[theorem]{Problem}
\newtheorem{example}[theorem]{Example}
\newtheorem{definition}[theorem]{Definition}
\newtheorem{remark}[theorem]{Remark}
\newtheorem{Fact}[theorem]{Fact}

\newcommand{\thm}[2]{\begin{theorem}\label{#1}{\sl #2}\end{theorem}}
\newcommand{\cor}[2]{\begin{corollary}\label{#1}{\sl #2}\end{corollary}}
\newcommand{\prop}[2]{\begin{proposition}\label{#1}{\sl #2}\end{proposition}}
\newcommand{\lem}[2]{\begin{lemma}\label{#1}{\sl #2}\end{lemma}}

\newcommand{\fact}[2]{\begin{Fact}\label{#1}{\sl #2}\end{Fact}}

\begin{document}
 
\maketitle

\begin{abstract} 
A function $F\colon\R^2\to\R$ is {\em sup-measurable\/}
if $F_f\colon\R\to\R$ given by $F_f(x)=F(x,f(x))$, $x\in\R$, 
is measurable for each measurable function $f\colon\R\to\R$. 
It is known that under different set theoretical assumptions, 
including CH, there are sup-measurable non-measurable 
functions, as well as their category analog. 
In this paper we will show
that the existence of category analog of 
sup-measurable non-measurable functions is independent of ZFC.
A similar result for the original measurable case 
is a subject of a work in prepartion by
Ros{\l}anowski and Shelah. 
\end{abstract}

\section{Introduction} 

Our terminology is standard and follows that from~\cite{BJ},
\cite{CiBook}, \cite{Ku}, or~\cite{ShPF}.

The study of sup-measurable 
functions comes from the theory of 
differential equations.
More precisely it comes from 
a question for which functions
$F\colon\R^2\to\R$
the Cauchy problem
\begin{equation} \label{Cy}
y'=F(x,y),\;\;\;y(x_0)=y_0
\end{equation}
has a (unique) {\em a.e.-solution\/} in the class
of locally absolutely continuous functions on $\real$
in a sense that 
$y(x_0)=y_0$ and $y'(x)=G(x,y(x))$ for 
almost all $x\in\real$. 
(For more on this motivation see \cite{Kh2} or \cite{BC}.)
It is not hard to find measurable functions which are not sup-measurable. 
(See \cite{Sr} or \cite[Cor.~1.4]{B}.)
Under the continuum hypothesis
CH or some weaker set-theoretical assumptions 
nonmeasurable sup-measurable functions were constructed 
in \cite{GL}, \cite{Kh1}, \cite{B}, and \cite{Kh2}. 
An independence from ZFC of the existence of such example 
is a subject of a work in prepartion by
Ros{\l}anowski and Shelah. 

A function $F\colon\R^2\to\R$ 
is a {\em category analog of 
sup-measurable 
function\/} (or {\em Baire sup-measurable}) 
provided 
$F_f\colon\R\to\R$ given by $F_f(x)=F(x,f(x))$, $x\in\R$, 
has the Baire property 
for each function $f\colon\R\to\R$ with the Baire property. 
Baire sup-measurable function without the Baire property
has been constructed under CH in~\cite{GG}.
(See also \cite{B} and \cite{BC}.) 
The main goal of this paper is to show that
the existence of such functions cannot be proved in ZFC.
For this we need the following easy fact.
(See \cite[Prop. 1.5]{B}.)

\prop{pr1}{
The following conditions are equivalent. 
\begin{description}
\item[(i)] There is a Baire sup-measurable function
$F\colon\R^2\to\R$ without the Baire property. 

\item[(ii)] There is a function
$F\colon\R^2\to\R$ without the Baire property
such that $F_f$ has the Baire property for every 
Borel function $f\colon\R\to\R$.

\item[(iii)] There is a set 
$A\subset\R^2$ without the Baire property 
such that the projection 
$\proj(A\cap f)=\{x\in\real\colon \la x,f(x)\ra\in A\}$ 
has the Baire property for each 
Borel function $f\colon\R\to\R$. 

\item[(iv)] There is a Baire sup-measurable function
$F\colon\R^2\to\{ 0,1\}$ without the Baire property.
\end{description}
}\

The equivalence of (i) and (ii) follows from the fact that
the function $F\colon\R^2\to\R$ is Baire sup-measurable
if and only if $F_f$ has a Baire property for every 
Borel function $f\colon\R\to\R$.
(It is also true that 
$F\colon\R^2\to\R$ is Baire sup-measurable
provided $F_f$ has a Baire property for every 
Baire class one function $f\colon\R\to\R$, and that 
$F\colon\R^2\to\R$ is sup-measurable
provided $F_f$ is measurable for every 
continuous function $f\colon\R\to\R$.
See for example \cite[Lem. 1 and Rem. 1]{BC}.)

The main theorem of the paper is the following.

\thm{thMain}{It is consistent with the set theory ZFC that 
for every $A\subset 2^\omega\times 2^\omega$ 
for which the sets $A$ and $A^c=(2^\omega\times 2^\omega)\setminus A$
are nowhere meager in 
$2^\omega\times 2^\omega$ 
there is a homeomorphism 
$f$ from $2^\omega$ onto $2^\omega$ such that
the set 
$\proj(A\cap f)$ does not have the  
Baire property in $2^\omega$.}

Before proving this theorem let us notice that
it implies easily the following corollary. 

\cor{cor1}{
The existence of Baire sup-measurable function
$F\colon\R^2\to\R$
without the Baire property is independent from 
the set theory ZFC. 
}

\proof Since under CH there are 
Baire sup-measurable functions
without the Baire property
it is enough to show that Theorem~\ref{thMain}
implies consistency with ZFC that there 
are no such functions. For this we will work
in the model from Theorem~\ref{thMain}.

So, take an arbitrary $A\subset\R^2$ without the Baire property.
By (iii) from Proposition~\ref{pr1} it
is enough to show there exists 
a Baire class one 
function $f\colon\R\to\R$ for which the set 
$\proj(A\cap f)$ does not have the Baire property. 

We will first show this under the additional assumption
that the sets $A$ and $\real^2\setminus A$
are nowhere meager in $\real^2$. 
But then the set $A_0=A\cap(\real\setminus\rational)^2$
and its complement are 
nowhere meager in $(\real\setminus\rational)^2$.
Moreover, since
$\real\setminus\rational$ is homeomorphic to
$2^\omega\setminus E$ for some countable set $E$
(the set of all eventually constant functions in $2^\omega$)
we can consider $A_0$ as a subset of 
$(2^\omega\setminus E)^2\subset 2^\omega\times 2^\omega$.
Then $A_0$ and its complement are 
still nowhere meager in $2^\omega\times 2^\omega$.
Therefore,  
there exists an autohomeomorphism 
$f$ of $2^\omega$
such that the set 
$\proj(A_0\cap f)=
\{x\in 2^\omega\setminus E\colon\la x,f(x)\ra\in A_0\}$  
does not have the  
Baire property in $2^\omega$.
Now, as before, $f\restriction(2^\omega\setminus E)$
can be considered as defined on 
$\real\setminus\rational$. 
So if $\bar f\colon\real\to\real$ is an extension
of $f\restriction(2^\omega\setminus E)$
(under such identification) to $\real$ 
as a constant on $\rational$  then $\bar f$ 
is Borel and the set 
$\proj(A_0\cap \bar f)$  
does not have the Baire property in $\real$.

Now, if $A$ is an arbitrary subset of 
$\R^2$ without the Baire property we can find 
non-empty open intervals 
$U$ and $W$ in $\real$ such that $A$ and $(U\times W)\setminus A$
are nowhere meager 
in $U\times W$. 
Since $U$ and $W$ are homeomorphic with $\real$ 
the above case implies the existence of
Borel function $f_0\colon U\to W$ such that
$\proj(A\cap f_0)$  
does not have the Baire property in $U$.
So any Borel extension $f\colon\real\to\real$ of $f_0$ 
works. 
\qed

\section{Reduction of the proof of Theorem~\ref{thMain}
to the main lemma}

The idea of the proof is quite simple. 
For every nowhere meager $A\subset 2^\omega\times 2^\omega$ 
for which $A^c=(2^\omega\times 2^\omega)\setminus A$
is also nowhere meager
we will find 
a natural ccc forcing notion $Q_A$ 
which adds the required homeomorphism 
$f$.
Then we will start with the constructible universe 
$V=L$ and 
iterate with finite support these 
notions of forcing in such a way that every
nowhere meager  
set $A^*\subset 2^\omega\times 2^\omega$,
with $(2^\omega\times 2^\omega)\setminus A^*$ nowhere meager, 
will be taken care of by some $Q_A$ 
at an appropriate step of iteration. 

There are two technical problems with carrying through this idea. 
First is that we cannot possible list in our iteration 
all nowhere meager  
subsets of 
$2^\omega\times 2^\omega$ with nowhere meager complements 
since the iteration can be of length at most 
continuum $\cont$ and there are $2^\cont$ many of such 
sets. This problem will be solved by defining our
iteration as 
$P_{\omega_2}=\la\la P_\alpha,\dot Q_\alpha\ra\colon\alpha<\omega_2\ra$ 
such that the generic extension $V[G]$ of $V$ with respect to 
$P_{\omega_2}$
will satisfy $2^\omega=2^{\omega_1}=\omega_2$ and
have the property that 
\begin{description}
\item{(m)} every non-Baire subset $A^*$ of $2^\omega$ contains 
           a non-Baire subset $A$ of cardinality $\omega_1$.
\end{description}
Thus in the iteration we will use only the forcing notions 
$Q_\alpha=Q_A$ 
for the sets $A$ of cardinality $\omega_1$, 
whose number is equal to $\omega_2$, the length of iteration.
Condition (m) will guarantee that this will give us enough
control of all nowhere meager  
subsets $A^*$ of $2^\omega\times 2^\omega$.

The second problem is that 
even if at some stage $\alpha<\omega_2$ of our iteration
we will add a homeomorphism 
$f$ appropriate
for a  given 
set $A\subset 2^\omega\times 2^\omega$,
that is such that
\[
V[G_\alpha]\models ``
\proj(A\cap f) \textrm{ is not Baire in $2^\omega$,''}  
\]
where $G_\alpha=G\cap P_\alpha$, 
then in general there is no guarantee that the set 
$\proj(A\cap f)$ will remain non-Baire at the final model $V[G]$.
The preservation of non-Baireness of each appropriate 
set $\proj(A\cap f)$ 
will be achieved by careful crafting our iteration following
a method known as 
the {\em oracle-cc\/} forcing iteration.

The theory of the oracle-cc forcings
is well described in \cite[Ch.~IV]{ShPF}
and here we will recall only the fragments 
that are relevant to our specific situation.
In particular if $\Gamma$ stands for the 
set of all limit ordinals 
less than $\omega_1$ then 
\begin{itemize}
\item an {\em $\omega_1$-oracle\/} is any sequence 
$\M=\la M_\delta\colon\delta\in\Gamma\ra$
where $M_\delta$ is a countable transitive model of ZFC$^-$
(i.e., 
ZFC without the power set axiom) such that 
$\delta\subset M_\delta$, 
$M_\delta\models ``\delta$\textrm{ is countable},''
and the set $\{\delta\in\Gamma\colon A\cap\delta\in M_\delta\}$
is stationary in $\omega_1$ for every $A\subset\omega_1$. 
\end{itemize}
The existence of an $\omega_1$-oracle is equivalent to 
the diamond principle $\diamondsuit$. 

We will also need the following fact which, for our purposes, can be
viewed as a definition of $\M$-cc property.  

\fact{f1}{If $P$ is a forcing notion of cardinality 
$\leq\omega_1$, 
$e\colon P\to\omega_1$ is one-to-one,
$\M=\la M_\delta\colon\delta\in\Gamma\ra$
is an $\omega_1$-oracle, and $C$ is a closed unbounded subset of
$\omega_1$
such that for every $\delta\in\Gamma\cap C$
\begin{quote}
$e^{-1}(E)$ is predense in $P$ for every set 
$E\in M_\delta\cap\P(\delta)$ for which 
$e^{-1}(E)$ is predense in $e^{-1}(\{\gamma\colon \gamma<\delta\})$
\end{quote}
then $P$ has the $\M$-cc property.
}
This follows immediately from the definition
of $\M$-cc property 
\cite[Definition~1.5, p. 119]{ShPF}
and \cite[Claim~1.4(3), p. 118]{ShPF}.


Our proof will relay on the following main lemma.

\lem{lemMain}{For every 
$A\subset 2^\omega\times 2^\omega$
for which $A$ and $A^c=(2^\omega\times 2^\omega)\setminus A$
are nowhere meager in $2^\omega\times 2^\omega$
and for every $\omega_1$-oracle $\M$
there exists an $\M$-cc forcing notion $Q_A$
of cardinality $\omega_1$ such that $Q_A$ forces
\begin{quote}
there exists 
an autohomeomorphism $f$ of 
$2^\omega$ such that the sets 
$\proj(f\cap A)$ and $\proj(f\setminus A)$ are nowhere 
meager in $2^\omega$.
\end{quote}
}

The proof of Lemma~\ref{lemMain} represents the core of our
argument and will be presented in the next section.
In the reminder of this section
we will sketch how Lemma~\ref{lemMain} 
implies Theorem~\ref{thMain}.
Since this follows the standard path,
as described in~\cite{ShPF},
the experts familiar with this treatment may 
proceed directly to the next section. 

Now, the 
iteration $P_{\omega_2}$ is defined 
by choosing by induction 
the sequence
$\la \la P_\alpha,\dot A_\alpha,\dot\M_\alpha,
\dot Q_\alpha, 
\dot f_\alpha
\ra\colon\alpha<\omega_2\ra$ 
such that for every $\alpha<\omega_2$
\begin{itemize}
\item $P_\alpha=\la \la P_\beta,\dot
      Q_\beta\ra\colon\beta<\alpha\ra$ is a finite support iteration,

\item $\dot A_\alpha$ is a $P_\alpha$-name 
and for every $\beta\leq\alpha$
\begin{itemize}
\item[] $P_\alpha\forces 
\mbox{``$\dot A_\beta$ and $(\dot A_\beta)^c$
are 
nowhere meager subsets of $2^\omega\times 2^\omega$,''}$
\end{itemize}

\item $\dot\M_\alpha$ is a $P_\alpha$-name such that
$P_\alpha$ forces 
\begin{itemize}
\item[] $\dot\M_\alpha$ is an $\omega_1$-oracle and 
for every $\beta<\alpha$
if $Q$ satisfies $\dot\M_\alpha$-cc then 
\[
\mbox{$Q\forces 
``\proj(\dot f_\beta \cap\dot A_\beta),
\proj(\dot f_\beta \setminus\dot A_\beta)\subset 2^\omega$ 
are nowhere meager 
in $2^\omega$,''}
\]
\end{itemize}

\item $\dot Q_\alpha$ is a $P_\alpha$-name 
for a forcing such that $P_\alpha$ forces
\begin{itemize}
\item[] $\dot Q_\alpha$ is an $\dot\M_\alpha$-cc forcing 
$Q_{\dot A_\alpha}$
from Lemma~\ref{lemMain},
\end{itemize}

\item $\dot f_\alpha$ is a $P_{\alpha+1}$-name 
for which $P_{\alpha+1}$ forces that 
\begin{itemize}
\item[] $\dot f_\alpha$ is a 
$\dot Q_\alpha$-name for the function $f$ from Lemma~\ref{lemMain}.
\end{itemize}
\end{itemize}
The existence of appropriate $\omega_1$-oracles
and the fact that each $P_\alpha$ obtained 
that way preserves nowhere meagerness
(i.e., 
non-meagerness of their traces on all basic open sets) 
of all the projections 
$\proj(\dot f_\beta \cap\dot A_\beta)$ and
$\proj(\dot f_\beta \setminus\dot A_\beta)$
for \mbox{$\beta<\alpha$} 
follows from Example~2.2 and results from section~3 
of \cite[ch.~IV]{ShPF}. Also, Claim~3.2 from
\cite[ch.~IV]{ShPF} implies that
all sets 
$\proj(\dot f_\alpha\cap\dot A_\alpha)$ and
$\proj(\dot f_\alpha\setminus\dot A_\alpha)$
remain 
nowhere meager in $2^\omega$ in the final model $V[G]$. 
Thus it is enough to 
ensure that each nowhere meager subset
$A^*$ of $2^\omega\times 2^\omega$
from $V[G]$ with nowhere meager complement 
contains an interpretation of 
some $\dot A_\alpha$. However, 
a choice of $\dot A_\alpha$'s which guarantee
this can be made with a help of the diamond principle 
$\diamondsuit_{\omega_2}$
and the fact that the set
\[
\{\alpha<\omega_2\colon \cf(\alpha)\neq\omega_1\mbox{ or }
V[G_\alpha]\models A\cap V[G_\alpha]\mbox{ is nowhere meager in }
2^\omega\times 2^\omega\}
\]
contains a closed unbounded set. (Compare~\cite[Claim~4.4, p. 130]{ShPF}.)

\section{Proof of Lemma~\ref{lemMain}}

Let $\K$ be the family of all sequences
$\bar h=\la h_\xi\colon\xi\in\Gamma\ra$ such that 
each $h_\xi$ is a function from
a countable set $D_\xi\subset 2^\omega$ onto
$R_\xi\subset 2^\omega$ and that
\[
D_\xi\cap D_\eta=R_\xi\cap R_\eta=\emptyset
\ \mbox{ for every distinct }\ \xi,\eta\in\Gamma.
\]
For each $\bar h\in\K$ we will define a forcing notion
$Q_{\bar h}$. 
Forcing $Q_A$ satisfying Lemma~\ref{lemMain}
will be chosen as $Q_{\bar h}$ for some $\bar h\in\K$.

So fix an $\bar h\in\K$. Then $Q_{\bar h}$ is defined
as the set of all triples $p=\la n,\pi,h\ra$ for which 

\begin{itemize}
\item[(A)] $h$ is a function from a finite subset $D$ of
$\bigcup_{\xi\in\Gamma}D_\xi$ into $2^\omega$;

\item[(B)] $n<\omega$ and $\pi$ is a permutation of $2^n$;

\item[(C)] $|D\cap D_\xi|\leq 1$ for every $\xi\in\Gamma$;

\item[(D)] if $x\in D\cap D_\xi$ then $h(x)=h_\xi(x)$ 
      and $h(x)\restriction n=\pi(x\restriction n)$.
\end{itemize}
Forcing $Q_{\bar h}$ is ordered as follows. 
Condition $p'=\la n',\pi',h'\ra$ is stronger than 
$p=\la n,\pi,h\ra$, $p'\leq p$, 
provided 
\begin{equation}\label{ord}
\mbox{$n\leq n'$, \ $h\subset h'$,\ \  and\ \ 
$\pi'(\eta)\restriction n=\pi(\eta\restriction n)$\ \ 
for every $\eta\in 2^{n'}$.}
\end{equation}

In the reminder of this paper we will write $[\eta]$
for the basic open neighborhood in $2^\omega$ generated by 
$\eta\in 2^{<\omega}$, that is, 
\[
[\eta]=\{x\in 2^\omega\colon \eta\subset x\}.
\]
Note that using this notation 
the second part of the condition (D) says that for every 
$x\in D$ and $\eta\in 2^n$
\begin{equation}\label{conD}
x\in [\eta]\ \  \mbox{ if and only if  }\ \  h(x)\in [\pi(\eta)].
\end{equation}
Also, if $n<\omega$ we will write $[\eta]\restriction 2^n$
for $\{x\restriction 2^n\colon x\in[\eta]\}$. 
Note that in this notation the part of (\ref{ord}) concerning 
permutations says that
$\pi'$ expends $\pi$ in a sense that $\pi'$ 
maps $[\zeta]\restriction 2^{n'}
$ onto 
$[\pi(\zeta)]\restriction 2^{n'}
$
for every $\zeta\in 2^{n}$.

In what follows we will use the following basic property of $Q_{\bar h}$.

\begin{itemize}
\item[($*$)] For every 
$q=\la n,\pi,h\ra\in Q_{\bar h}$ and $m<\omega$
there exist an $n'\geq m$ and a permutation $\pi'$ of $2^{n'}$
such that 
$q'=\la n',\pi',h\ra
\in Q_{\bar h}$ and $q'$ extends $q$. 
\end{itemize}

The choice of such $n'$ and $\pi'$ is easy. 
First pick $n'\geq \max\{m,n\}$ such that
$x\restriction n'\neq y\restriction n'$
for every different $x$ and $y$ from either domain $D$ 
or range $R=h[D]$ of $h$.
This implies that that for every $\zeta\in 2^{n}$
the set 
$D_\zeta=
\{x\restriction n'\colon x\in D\cap [\zeta]\}\subset [\zeta]\restriction 2^{n'}$ 
has the same cardinality that $D\cap [\zeta]$
and $H_\zeta=
\{x\restriction n'\colon x\in h[D]\cap [\pi(\zeta)]\}
\subset [\pi(\zeta)]\restriction 2^{n'}$ 
has the same cardinality that $h[D]\cap [\pi(\zeta)]$.
Since, by (\ref{conD}), we have also
$|D\cap [\zeta]|=|h[D]\cap [\pi(\zeta)]|$
we see that $|D_\zeta|=|H_\zeta|$. 
Define $\pi'$ on $D_\zeta$ by 
$\pi'(x\restriction n')=h(x)\restriction n'$
for every $x\in D_\zeta$.  Then $\pi'$ is a bijection 
from $D_\zeta$ onto $H_\zeta$ 
and this definition ensures that an appropriate part of
the condition (D) for $h$ and $\pi'$ is
satisfied.  Also, if for each $\zeta\in 2^{n}$
we extend $\pi'$ onto $[\zeta]\restriction 2^{n'}$
as a bijection from $([\zeta]\restriction 2^{n'})\setminus D_\zeta$
onto $([\pi(\zeta)]\restriction 2^{n'})\setminus H_\zeta$,
then the condition (\ref{ord}) will be satisfied. 
Thus such defined $q'=\la n',\pi',h\ra$ belongs to 
$Q_{\bar h}$ and extends $q$.

Next note that forcing $Q_{\bar h}$ 
has the following properties needed to prove Lemma~\ref{lemMain}. 
In what follows we will 
consider $2^\omega$ with the standard
distance: 
$$
d(r_0,r_1)=2^{-\min\{n<\omega\colon r_0(n)\neq r_1(n)\}}
$$
for different $r_0,r_1\in 2^\omega$.

\fact{fact2}{Let $\bar h=\la h_\xi\colon\xi\in\Gamma\ra\in\K$ and 
$f=\bigcup\{h\colon \la n,\pi,h\ra\in H\}$, where 
$H$ is a $V$-generic filter over $Q_{\bar h}$.
Then $f$ is a uniformly continuous one-to-one function from a subset
$D$ of $2^\omega$ into~$2^\omega$.
Moreover, if for every $\xi\in\Gamma$ the graph of $h_\xi$
is dense in $2^\omega\times 2^\omega$  then $D$ and $f[D]$ are dense in
$2^\omega$ and $f$ can be uniquely extended to an autohomeomorphism 
$\tilde f$ of $2^\omega$. 
}

\proof Clearly $f$ is a one-to-one function from a subset
$D$ of $2^\omega$ into~$2^\omega$.
To see that it is uniformly continuous 
choose an $\ep>0$.
We will find $\delta>0$ such that 
$r_0,r_1\in D$ and $d(r_0,r_1)<\delta$ imply $d(f(r_0),f(r_1))<\ep$. 
For this note that, by ($*$), the set
\[
S=\{q=\la n,\pi,h\ra\in Q_{\bar h}\colon 2^{-n}<\ep\}
\]
is dense in $Q_{\bar h}$. So take $q=\la n,\pi,h\ra\in H\cap S$
and put $\delta=2^{-n}$.
We claim that this $\delta$ works.

Indeed, take $r_0,r_1\in D$ with $d(r_0,r_1)<\delta$.
Then there is $q'=\la n',\pi',h'\ra\in H$ 
stronger than $q$ 
such that $r_0$ and $r_1$ are in the domain of $h'$.
Therefore, $n\leq n'$ and for $j<2$
\[
f(r_j)\restriction n
=h'(r_j)\restriction n
=(h'(r_j)\restriction n')\restriction n
=\pi'(r_j\restriction n')\restriction n
=\pi(r_j\restriction n)
\]
by conditions (D) and (\ref{ord}).
Since $d(r_0,r_1)<\delta=2^{-n}$ 
implies $r_0\restriction n=r_1\restriction n$
we obtain
\[
f(r_0)\restriction n
=\pi(r_0\restriction n)
=\pi(r_1\restriction n)
=f(r_1)\restriction n
\]
that is, 
$d(f(r_0),f(r_1))<2^{-n}<\ep$. So $f$ is uniformly continuous. 

Essentially the same argument (with the same values of $\ep$ and $\delta$)
shows that $f^{-1}\colon f[D]\to D$ is uniformly continuous. 
Thus, if $\tilde f$ is the unique continuous extension of 
$f$ into $\cl(D)$ then $\tilde f$ is a homeomorphism from
$\cl(D)$ onto $\cl f[D]$.

To finish the argument assume that all functions $h_\xi$ have dense
graphs, take an $\eta\in 2^{m}$ for some $m<\omega$,
and notice that the set
\[
S_\eta=\{q=\la n,\pi,h\ra\in Q_{\bar h}\colon 
\mbox{ the domain $D'$ of $h$ intersects }[\eta]\}
\]
is dense in $Q_{\bar h}$. Indeed, if 
$q=\la n,\pi,h\ra\in Q_{\bar h}$
then, by ($*$), strengthening $q$ if necessary, 
we can assume that $m\leq n$.
Then, refining $\eta$ if necessary, we can also assume that $m=n$,
that is, that $\eta$ is in the domain of $\pi$.
Now, if $[\eta]$ intersects the domain of 
$h$ then already $q$ belongs to $S_\eta$.
Otherwise take $\xi\in\Gamma$ with $D'\cap D_\xi=\emptyset$
and pick $\la x,h_\xi(x)\ra\in [\eta]\times[\pi(\eta)]$,
which exists by the density of the graph of $h_\xi$. 
Then $\la n,\pi,h\cup\{\la x,h_\xi(x)\ra\}\ra$ belongs to $S_\eta$ 
and extends $q$. 

This shows that $D\cap [\eta]\neq\emptyset$ for every 
$\eta\in 2^{<\omega}$, that is, $D$ is dense in $2^\omega$. 

The similar argument shows that for every 
$\eta\in 2^{<\omega}$ the set
\[
S^\eta=\{q=\la n,\pi,h\ra\in Q_{\bar h}\colon 
\mbox{ the range of $h$ intersects }[\eta]\}
\]
is dense in $Q_{\bar h}$, which implies that  
$h[D]$ is dense in $2^\omega$. 
Thus $\tilde f$ is a homeomorphism from $\cl(D)=2^\omega$
onto $\cl h[D]=2^\omega$.
\qed

Now take $A\subset 2^\omega\times 2^\omega$
for which $A$ and $A^c=(2^\omega\times 2^\omega)\setminus A$
are nowhere meager in $2^\omega\times 2^\omega$
and fix an $\omega_1$-oracle 
$\M=\la M_\delta\colon\delta\in\Gamma\ra$. 
By Fact~\ref{fact2} in order to prove Lemma~\ref{lemMain}
it is enough to find an 
$\bar h=\la h_\xi\colon\xi\in\Gamma\ra\in\K$ such that
\begin{equation}\label{eqN2}
\mbox{$Q_A=Q_{\bar h}$ is $\M$-cc} 
\end{equation}
and $Q_{\bar h}$ forces that, in $V[H]$, 
\begin{equation}\label{eqN3}
\mbox{the sets 
$\proj(f\cap A)$ and $\proj(f\setminus A)$ are nowhere 
meager in $2^\omega$.} 
\end{equation}
(In (\ref{eqN3}) function $f$ is defined as in Fact~\ref{fact2}.)

To define $\bar h$ we will construct 
a sequence
$\la \la x_\alpha,y_\alpha\ra\in 2^\omega\times 2^\omega\colon
\alpha<\omega_1\ra$ 
aiming for 
$h_\xi=\{\la x_{\xi+n},y_{\xi+n}\ra\colon n<\omega\}$,
where $\xi\in\Gamma$.

Let $\{\la\eta_n,\zeta_n\ra\colon n<\omega\}$
be an enumeration of
$2^{<\omega}\times 2^{<\omega}$. 
Points $\la x_{\xi+n},y_{\xi+n}\ra$
are chosen inductively in such a way that
\begin{description}
\item{(i)} $\la x_{\xi+n},y_{\xi+n}\ra$ 
      is a Cohen generic number over 
      $M_\delta[\la \la x_\alpha,y_\alpha\ra\colon \alpha<\xi+n\ra]$
      for every $\delta\leq\xi$, $\delta\in\Gamma$, 
      that is, $\la x_{\xi+n},y_{\xi+n}\ra$ is outside of all 
      meager subsets of $2^\omega\times 2^\omega$ which are
      coded in 
      $M_\delta[\la \la x_\alpha,y_\alpha\ra\colon \alpha<\xi+n\ra]$;
\item{(ii)} $\la x_{\xi+n},y_{\xi+n}\ra\in A$ if $n$ is even, and
      $\la x_{\xi+n},y_{\xi+n}\ra\in A^c$ otherwise. 
\item{(iii)} $\la x_{\xi+n},y_{\xi+n}\ra\in [\eta_n]\times[\zeta_n]$. 
\end{description}
The choice of $\la x_{\xi+n},y_{\xi+n}\ra$
is possible since both sets $A$ and $A^c$
are nowhere meager, and we consider
each time only countably many meager sets. 
Condition (iii) guarantees that the graph of each of $h_\xi$
will be dense in $2^\omega\times 2^\omega$. 
Note that if $\Gamma\ni\delta\leq \alpha_0<\cdots<\alpha_{k-1}$,
where $k<\omega$, then (by the product lemma in $M_\delta$)
\begin{equation}\label{PLem}
\mbox{$\la\la x_{\alpha_i},y_{\alpha_i}\ra\colon i<k\ra$
is an $M_\delta$-generic Cohen number in $\left(2^\omega\right)^k$.}
\end{equation}

For $q=\la n,\pi,h\ra\in Q_{\bar h}$
define
$$
\hat q=\bigcup_{\la\xi,\zeta\ra\in\pi}[\xi]\times[\zeta].
$$
Clearly $\hat q$ is an open subset of $2^\omega\times 2^\omega$
and condition (\ref{ord}) implies that 
for every $q,r\in Q_{\bar h}$ with $r=\la n',\pi',h'\ra$
\begin{equation}\label{BBBB}
\mbox{if $q\leq r$ then $\hat q\subset\hat r$
and $\hat q\cap([\xi]\times[\zeta])\neq\emptyset$ 
for every $\la\xi,\zeta\ra\in\pi'$.}
\end{equation}
Also for $\delta\in\Gamma$ let
$(Q_{\bar h})^\delta=\left\{\la n,\pi,h\ra\in Q_{\bar h}
\colon h\subset\bigcup_{\zeta<\delta}h_\zeta\right\}$.
To prove (\ref{eqN2}) and (\ref{eqN3}) we will use also the following fact.

\fact{FACT}{Let $\delta\in\Gamma$ be such that 
$M_\delta$ contains 
$(Q_{\bar h})^\delta$ and let $E\in M_\delta$ be, in $M_\delta$,
a predense subset of $(Q_{\bar h})^\delta$. Then
for every $k<\omega$ and
$p=\la n,\pi,h\ra\in (Q_{\bar h})^\delta$ 
the set 
\begin{equation}\label{CCC}
B_{p}^k=\bigcup\left\{
(\hat q)^k\colon 
\mbox{ $q$ extends $p$ and some $q_0\in E$}\right\}
\end{equation}
is dense in $(\hat p)^k\subset\left(2^\omega\times 2^\omega\right)^k$.
}

\proof 
By way of contradiction assume that $B_{p}^k$ is not dense in $(\hat p)^k$.
Then there are $t<\omega$ and 
$\xi_0,\zeta_0,\ldots,\xi_{k-1},\zeta_{k-1}\in 2^t$
such that $P=\prod_{i<k}([\xi_i]\times[\zeta_i])\subset (\hat p)^k$
is disjoint with $B_{p}^k$. 
Increasing $t$ and refining
$\xi_i$'s and $\zeta_j$'s, if necessary, we may assume that $t\geq n$, 
all $\xi_i$'s and $\zeta_j$'s are different, 
$\bigcup_{i<k}[\xi_i]$
is disjoint from the domain $D$ of $h$,
and $h[D]\cap\bigcup_{i<k}[\zeta_i]=\emptyset$. 
We can also assume that $x\restriction t\neq y\restriction t$ for every different
$x$ and $y$ from $D$ and from $h[D]$. Now, refining slightly the argument for  
($*$) we can find 
$r=\la t,\pi',h\ra\in (Q_{\bar h})^\delta$ extending $p$ such that
$\pi'(\xi_i)=\zeta_i$ for every $i<k$.
(Note that $P\subset (\hat p)^k$.)
We will obtain a contradiction with predensity of 
$E$ in $(Q_{\bar h})^\delta$ by showing that 
$r$ is incompatible with every element of $E$.

Indeed if $q$ were an extension of $r\leq p$ and an element $q_0$ of
$E$ then we would have 
$(\hat q)^k\subset B^k_p$.
But then, by (\ref{BBBB}) and the fact that $\la \xi_i,\zeta_i\ra\in\pi'$
for $i<k$,
we would also have $(\hat q)^k\cap P\neq\e$, 
contradicting $P\cap B^m_p=\e$. 
This finishes the proof of Fact~\ref{FACT}. \qed

Now we are ready to prove (\ref{eqN2}), that is, that 
$Q_{\bar h}$ is $\M$-cc. So, fix a bijection
$e\colon Q_{\bar h}\to\omega_1$ and let
\[
C=\left\{\delta\in\Gamma\colon (Q_{\bar h})^\delta=e^{-1}(\delta)\right\}.
\]
Then $C$ is a closed unbounded subset of $\omega_1$. 
Take a $\delta\in C$ for which $(Q_{\bar h})^\delta=e^{-1}(\delta)\in M_\delta$
and fix an $E\subset \delta$, $E\in M_\delta$, for which 
$e^{-1}(E)$ is predense in 
$(Q_{\bar h})^\delta$. 
By Fact~\ref{f1} it is enough to show that
\begin{equation*}\label{AAA}
\mbox{$e^{-1}(E)$ is predense in $Q_{\bar h}$.}
\end{equation*}

Take
$p_0=\la n,\pi,h_0\ra$ from $Q_{\bar h}$,
let $h=h_0\restriction\bigcup_{\eta<\delta}D_\eta$ and
$h_1=h_0\setminus h$, 
and notice that the condition 
$p=\la n,\pi,h\ra$ 
belongs to $(Q_{\bar h})^\delta$.
Assume that $h_1=\{\la x_i,y_i\ra\colon i<k\}$.
Since $s(h_1)=\la\la x_i,y_i\ra\colon i<k\ra\in (\hat p)^k$ and, by
Fact~\ref{FACT},  
$B_{p}^k$ is dense in $(\hat p)^k$
condition (\ref{PLem}) implies that 
$s(h_1)\in B_{p}^k$.
So there is
$q=\la n_0,\pi_0,g\ra\in(Q_{\bar h})^\delta$ 
extending $p$ and some $q_0\in e^{-1}(E)$ for which
$s(h_1)\in \hat q^m$.
But then $p'=\la n_0,\pi_0,g\cup h_1\ra$
belongs to $Q_{\bar h}$ and extends 
$q$.
This finishes the proof of (\ref{AAA}).

The proof of (\ref{eqN3}) is similar. We will prove only that
$\proj(f\setminus A)=\proj(f\cap A^c)$ is nowhere meager in
$2^\omega$, the argument for $\proj(f\cap A)$ being essentially the
same. 

By way of contradiction assume that $\proj(f\setminus A)$ is not
nowhere meager in $2^\omega$. So there is an $\eta\in 2^{<\omega}$
such that $\proj(f\setminus A)$ is meager in $[\eta]$. 
Thus, there is a sequence $\la \dot U_m\colon m<\omega\ra$
of $Q_{\bar h}$-names 
for which $Q_{\bar h}$ forces 
\[
\mbox{each $\dot U_m$ is an open dense subset of $[\eta]$ and 
$\proj(f\cap A)\cap\bigcap_{m<\omega}\dot U_m=\e$.
}
\]
Moreover, since $Q_{\bar h}$ is ccc (as every $\M$-cc forcing is ccc),
we can also assume that there exists a $\delta_0\in\Gamma$
such that 
$\la \dot U_m\colon m<\omega\ra\in M_{\delta_0}$.

Now, by the definition of $\omega_1$-oracle, the set 
\[
B_0=\{\delta\in\Gamma\colon 
\la \dot U_m\colon m<\omega\ra\in M_\delta\ \&\ 
M_\delta\cap Q_{\bar h}=(Q_{\bar h})^{\delta}\in M_\delta\}
\]
is stationary in $\omega_1$. Thus, we can find, in $V[H]$,
a $\delta\in B_0$ and an odd $j<\omega$ such that 
$x_{\delta+j}\in[\eta]\cap \proj(f)$.
Recall that $\la x_{\delta+j},y_{\delta+j}\ra \in A^c$ for odd $j$'s.
Therefore 
$x_{\delta+j}\in \proj(f\setminus A)$. Let $p_0=\la
n_0,\pi_0,h_0\ra\in H$ be such that 
$p_0\forces``x_{\delta+j}\in [\eta]\cap \proj(f).$''
We can assume that $x_{\delta+j}$ belongs to the domain of $h_0$. 
We will show that
\[
p_0\forces x_{\delta+j}\in \bigcap_{m<\omega}\dot U_m,
\]
which will finish the proof. 

So, assume that this is not the case. Then there exist
an $i<\omega$ and 
$p_1=\la n,\pi,h_1\ra\in Q_{\bar h}$ stronger than $p_0$ 
such that $p_1\forces``x_{\delta+j}\notin\dot U_i$.'' 
Let $h=h_1\restriction\bigcup_{\eta<\delta}D_\eta$ and
$h_1\setminus h=\{\la a_t,b_t\ra\colon t<m\}$.
Notice that 
$p=\la n,\pi,h\ra$ 
belongs to $(Q_{\bar h})^\delta$.
We can also assume that 
$\la x_{\delta+j},y_{\delta+j}\ra=\la a_0,b_0\ra$.

Now consider the set
\[
Z=\left\{q=\la n^q,\pi^q,h^q\ra\in(Q_{\bar h})^\delta\colon 
q\leq p\ \&\ 
\left(\exists \la \xi,\zeta\ra\in\pi^q\right) 
q\forces ``[\xi]\subset\dot U_i\mbox{''}\right\}
\]
and note that it belongs to $M_\delta$. 
Notice also that it is dense below $p$ in $(Q_{\bar h})^\delta$.
This is the case since any persistent information 
on $\dot U_i$ depends only on conditions from $(Q_{\bar h})^\delta$.
Thus, by Fact~\ref{FACT}, the set
\[
B_{p}^m=\bigcup\left\{
(\hat q)^m\colon 
\mbox{$q$ extends some $q_0\in Z$}\right\}\in M_\delta
\]
is dense in $(\hat p)^m$. So, by (\ref{PLem}), 
$\la\la a_t,b_t\ra\colon t<m\ra\in B_{p}^m$
since $\la\la a_t,b_t\ra\colon t<m\ra$
belongs to $(\hat p_1)^m=(\hat p)^m$.
But this means that there exist 
$q=\la n^q,\pi^q,h^q\ra\in Z$ and $\la \xi,\zeta\ra\in\pi^q$
with $q\forces ``[\xi]\subset\dot U_i\mbox{''}$
for which $\la\la a_t,b_t\ra\colon t<m\ra\in (\hat q)^m$.
A slight modification of the proof of Fact~\ref{FACT}
lets us also to choose $q$ such that 
$x_{\delta+j}=a_0\in[\xi]$.
But then $p_2=\la n^q,\pi^q,h^q\cup\{\la a_t,b_t\ra\colon t<m\}\ra$
belongs to $Q_{\bar h}$ and extends $p_1$. So,
$p_2$ forces that $x_{\delta+j}=a_0\in[\xi]\subset \dot U_i$
contradicting our assumption that 
$p_1\forces``x_{\delta+j}\notin\dot U_i$.'' 

This finishes the proof of (\ref{eqN3}) and 
and of Lemma~\ref{lemMain}.

\end{document}